\documentclass{article}

\newtheorem{thm}{Theorem}[section]

\newtheorem{lem}[thm]{Lemma}

\begin{document}

\title{{\bf An Improved Error Bound for Gaussian Interpolation}}         
\author{{\bf Lin-Tian Luh} \\Department of Mathematics, Providence University\\ Shalu, Taichung\\ email:ltluh@pu.edu.tw \\ phone:(04)26328001 ext. 15126 \\ fax:(04)26324653 }        
\date{\today}          
\maketitle
{\bf Abstract.} It's well known that there is a so-called exponential-type error bound for Gaussian interpolation which is the most powerful error bound hitherto. It's of the form $|f(x)-s(x)|\leq c_{1}(c_{2}d)^{\frac{c_{3}}{d}}\| f\| _{h}$ where $f$ and $s$ are the interpolated and interpolating functions respectively, $c_{1},c_{2},c_{3}$ are positive constants, $d$ is the fill-distance which roughly speaking measures the spacing of the data points, and $\| f\| _{h}$ is the $h$-norm of $f$ where $h$ is the Gaussian function. The error bound is suitable for $x\in R^{n},\ n\geq 1$, and gets small rapidly as $d\rightarrow 0$. The drawback is that the crucial constants $c_{2}$ and $c_{3}$ get worse rapidly as $n$ increases in the sense $c_{2}\rightarrow \infty$ and $c_{3}\rightarrow 0$ as $n\rightarrow \infty$. In this paper we raise an error bound of the form $|f(x)-s(x)|\leq c_{1}'(c_{2}'d)^{\frac{c_{3}'}{d}}\sqrt{d}\| f\| _{h}$, where $c_{2}'$ and $c_{3}'$ are independent of the dimension $n$. Moreover, $c_{2}'<<c_{2},\ c_{3}<<c_{3}'$, and $c_{1}'$ is only slightly different from $c_{1}$. What's important is that all constants $c_{1}',c_{2}'$ and $c_{3}'$ can be computed without slight difficulty.\\
\\
{\bf AMS classification}:41A05,41A25,41A30,41A63,65D10.\\
\\
{\bf keywords}:radial basis function, interpolation, error bound, Gaussian.
\section{Introduction}       
First, let $h$ be a continuous radial function on $R^{n}$ which is conditionally positive definite of order $m$. Given data points $(x_{j},f_{j}),j=1,\ldots , N$, where $X=\{ x_{1}, \ldots ,x_{N} \} $ is a subset of points in $R^{n}$ and the $f_{j}'s$ are real or complex numbers, the so-called $h$ spline interpolant of these data points is the function $s$ defined by 
\begin{equation}
  s(x)=p(x)+\sum_{j=1}^{N}c_{j}h(x-x_{j}),
\end{equation}
where $p(x)$ is a polynomial in $ P^{n}_{m-1}$ and $c_{j}'s$ are chosen so that
\begin{equation}
  \sum_{j=1}^{N}c_{j}q(x_{j})=0
\end{equation}
for all polynomials $q$ in $ P^{n}_{m-1}$ and 
\begin{equation}
  p(x_{i})+\sum_{j=1}^{N}c_{j}h(x_{i}-x_{j})=f_{i},\ \ i=1,\ldots ,N.
\end{equation}
Here $ P^{n}_{m-1}$ denotes the class of those polynomials of $R^{n}$ of degree $\leq m-1$.

A famous property is that the system of equations (2) and (3) has a unique solution whenever $X$ is a determining set for $ P^{n}_{m-1}$ and $h$ is strictly conditionally positive definite. More details can be seen in \cite{MN1}. Therefore in our case the interpolant $s(x)$ is well defined.

We clarify that $X$ is said to be a determining set for $ P^{n}_{m-1}$ if $X$ does not lie on the zero set of any nontrivial polynomial in $ P^{n}_{m-1}$.

In this paper the function $h$ is defined by
\begin{equation}
  h(x):=e^{-\beta |x|},\ \beta>0.
\end{equation} 
This is the so-called Gaussian function.

In \cite{MN3} Madych and Nelson raise the famous exponential-type error bound for the scattered data interpolation of Gaussian function, as mentioned in the abstract. The computation of the constants $c_{1},c_{2}$ and $c_{3}$ can be found in \cite{Lu3}. This error bound is very powerful. However, too many data points may lead to a large condition number when solving the linear system (2) and (3). There will be a significant improvement of the ill-conditioning if a better error bound can be obtained so that a satisfactory error estimate is reached before too many data points are involved. This is what we are pursuing in this paper.
\subsection{Polynomials and Simplices} 
In this paper we will use $P_{l}^{n}$ to denote the space of polynomials of degree $l$ in $n$ variables. It's well known that it has dimension $dimP_{l}^{n}=\left( \begin{array}{c}
                                                                        n+l \\ n   
                                                                      \end{array} \right)$. We will denote $dimP_{l}^{n}$ by $N$ in this section. Let $E\subseteq R^{n}$ be compact. The interpolation theory tells us that if $x_{1},\ldots ,x_{N} \in E$ and do not lie on the zero set of any nontrivial $q\in P_{l}^{n}$, there exists Lagrange polynomials $l_{i},\ i=1,\ldots ,N$, of degree $l$ defined by $l_{i}(x_{j})=\delta_{ij},\ 1\leq i,j\leq N$, such that for any $f\in C(E),\ (\Pi_{l}f)(x):=\sum_{i=1}^{N}f(x_{i})l_{i}(x)$ is its interpolating polynomial. It's easily seen that $\Pi_{l}(p)=p$ for all $p\in P_{l}^{n}$ and hence the mapping $\Pi_{l}:C(E)\rightarrow P_{l}^{n}$ is a projection. Let
$$\| \Pi_{l}\| :=\max _{x\in E}\sum_{i=1}^{N}|l_{i}(x)|.$$
A famous property says that for any $p\in P_{l}^{n}$,

$$\| p\| _{\infty}\leq \| \Pi_{l}\| \max _{1\leq i\leq N}|p(x_{i})|$$
The compact set $E$ discussed in this paper will be mainly an $n$-dimensional simplex $T_{n}$ whose definition can be found in \cite{Fl}.

It will be convenient to adopt {\bf barycentric coordinates} when discussing points in a simplex. Suppose $v_{1},\ldots ,v_{n+1}$ are the vertices of $T_{n}$. Then any $x\in T_{n}$ can be written as a convex combination of the vertices:
$$x=\sum_{i=1}^{n+1}c_{i}v_{i}$$
 where $\sum_{i=1}^{n+1}c_{i}=1$ and $c_{i}\geq 0$ for all $i$. The barycentric coordinate of $x$ is then $(c_{1},\ldots ,c_{n+1})$. Let's define ``{\bf equally spaced}'' points of {\bf degree} $l$ to be those points whose barycentric coordinates are of the form
$$(k_{1}/l,k_{2}/l,\ldots ,k_{n+1}/l),\ k_{i}\ nonnegative \ integers\ with\ \sum_{i=1}^{n+1}k_{i}=l.$$
Obviously the number of such points in $T_{n}$ is exactly $N=dimP_{l}^{n}$. Moreover, by \cite{Bo}, we know that equally spaced points form a determiming set for $P_{l}^{n}$.

The following lemma cited from \cite{Bo} will be needed.
\begin{lem}
  For the above equally spaced points $\| \Pi_{l}\| \leq \left( \begin{array}{c}
                                                                  2l-1 \\ l
                                                                \end{array} \right) $. Moreover, as $n\rightarrow \infty ,\ \| \Pi_{l}\| \rightarrow \left( \begin{array}{c}
                                                             2l-1 \\ l
                                                           \end{array} \right) $.
\end{lem}
Now, we are going to prove a lemma which plays a crucial role in our construction of the error bound.
\begin{lem}
  Let $Q\subseteq R^{n}$ be an n simplex in $R^{n}$ and $Y$ be the set of equally spaced points of degree $l$ in $Q$. Then, for any point $x$ in $Q$, there is a measure $\sigma$ supported on $Y$ such that 
$$\int p(y)d\sigma(y)=p(x)$$
for all $p$ in $P_{l}^{n}$, and 
$$\int d|\sigma |(y)\leq \left( \begin{array}{c}
                                2l-1 \\ l
                              \end{array} \right).$$
\end{lem}
{\bf Proof}. Let $Y=\{ y_{1},\ldots , y_{N}\} $ be the set of equally spaced points of degree $l$ in $Q$. Denote $P_{l}^{n}$ by $V$. For any $x\in Q$, let $\delta_{x}$ be the point-evaluation functional. Define $T:V\rightarrow T(V)\subseteq R^{N}$ by $T(v)=(\delta_{y_{i}}(v))_{y_{i}\in Y}$. Then $T$ is injective. Define $\tilde{\psi}$ on $T(V)$ by $\tilde{\psi}(w)=\delta_{x}(T^{-1}w)$. By the Hahn-Banach theorem, $\tilde{\psi}$ has a norm-preserving extension $\tilde{\psi}_{ext}$ to $R^{N}$. By the Riesz representation theorem, each linear functional on $R^{N}$ can be represented by the inner product with a fixed vector. Thus, there exists $z\in R^{N}$ with 
$$\tilde{\psi}_{ext}(w)=\sum_{j=1}^{N}z_{j}w_{j}$$
and $\| z\| _{(R^{N})^{*}}=\| \tilde{\psi}_{ext}\| $. If we adopt the $l_{\infty}$-norm on $R^{N}$, the dual norm will be the $l_{1}$-norm. Thus $\| z\| _{(R^{N}))^{*}}=\| z\| _{1}=\| \tilde{\psi}_{ext}\| =\| \tilde{\psi}\| =\| \delta_{x}T^{-1}\|$.

Now, for any $p\in V$, by setting $w=T(p)$, we have 
$$\delta_{x}(p)=\delta_{x}(T^{-1}w)=\tilde{\psi}(w)=\tilde{\psi}_{ext}(w)=\sum_{j=1}^{N}z_{j}w_{j}=\sum_{j=1}^{N}z_{j}\delta_{y_{j}}(p).$$
This gives 
\begin{equation}
  p(x)=\sum_{j=1}^{N}z_{j}p(y_{j})
\end{equation}
where $|z_{1}|+\cdots +|z_{N}|=\| \delta_{x}T^{-1}\|$.

Note that
\begin{eqnarray*}
  \| \delta_{x}T^{-1}\| & = & \sup _{\begin{array}{c}
                                       w\in T(V) \\ w\neq 0
                                     \end{array} }\frac{\| \delta_{x}T^{-1}(w)\| }{\| w\| _{R^{N}}}\\
                        & = & \sup _{\begin{array}{c}
                                       w\in T(V) \\ w\neq 0
                                     \end{array}}\frac{|\delta_{x}p|}{\| T(p)\| _{R^{N}}}\\
                        & \leq & \sup _{\begin{array}{c}
                                          p\in V \\ p\neq 0
                                        \end{array}}\frac{|p(x)|}{\max _{j=1,\ldots ,N}|p(y_{j})|}\\
                        & \leq & \sup _{\begin{array}{c}
                                          p\in V \\ p\neq 0
                                        \end{array}}\frac{\| \Pi_{l}\| \max _{j=1,\ldots ,N}|p(y_{j})|}{\max _{j=1,\ldots ,N}|p(y_{j})|}\\
                        & = & \| \Pi_{l}\| \\
                        & \leq & \left( \begin{array}{c}
                                          2l-1 \\ l
                                        \end{array}      \right) .
                                                                                                  \end{eqnarray*}
Therefore $|z_{1}|+\cdots +|z_{N}|\leq \left( \begin{array}{c}
                                               2l-1 \\ l
                                             \end{array} \right)$ and our lemma follows immediately by (5). \hspace{2cm}  $\sharp$
\subsection{Radial Functions and Borel Measures}
Before moving on to our main result, some background for interpolation is necessary. First, the space of complex-valued functions on $R^{n}$ that are compactly supported and infinitely differentiable is denoted by ${\cal D}$. The Fourier transform of a function $\phi$ in ${\cal D}$ is
$$\hat{\phi}(\xi)=\int e^{-i<x,\xi >}\phi(x)dx.$$ 
Then a crucial lemma introduced in \cite{GV} but modified by Madych and Nelson in \cite{MN2} says that for any continuous conditionally positive definite function $h$ of order $m$, the Fourier transform of $h$ uniquely determines a positive Borel measure $\mu$ on $R^{n}\sim \{ 0\} $ and constants $a_{\gamma}, |\gamma|=2m$ as follows: For all $\psi \in {\cal D}$
\begin{eqnarray}
  \int h(x)\psi(x)dx & = & \int \left\{ \hat{\psi}(\xi)-\hat{\chi}(\xi)\sum_{|\gamma|<2m}D^{\gamma}\hat{\psi}(0)\frac{\xi^{\gamma}}{\gamma !}\right\} d\mu(\xi) \nonumber \\
                     &   & + \sum_{|\gamma|\leq 2m}D^{\gamma}\hat{\psi}(0)\frac{a_{\gamma}}{\gamma!}, 
\end{eqnarray}
where for every choice of complex numbers $c_{\alpha},|\alpha|=m$,
$$\sum_{|\alpha|=m}\sum_{|\beta|=m}a_{\alpha+\beta}c_{\alpha}\overline{c_{\beta}}\geq 0.$$
Here $\chi$ is a function in ${\cal D}$ such that $1-\hat{\chi}(\xi)$ has a zero of order $2m+1$ at $\xi=0$; both of the integrals
$$\int_{0<|\xi|<1}|\xi |^{2m}d\mu(\xi),\ \ \int_{|\xi|\geq 1}d\mu(\xi)$$
are finite. The choice of $\chi$ affects the value of the coefficients $a_{\gamma}$ for $|\gamma|<2m$.\\
\section{Main Result}
  Before showing our main result, we need some lemmas. First, recall the famous formula of Stirling.\\
\\
{\bf Stirling's Formula}: $n!\sim \sqrt{2\pi n}(\frac{n}{e})^{n}$.\\
\\
The approximation is very reliable even for small $n$. For example, when $n=10$, the relative error is only $0.83\%$ . The larger $n$ is, the better the approximation is. For further details, we refer the reader to \cite{GG} and \cite{GKP}.
\begin{lem}
  Let $\rho_{1}=\frac{1}{e}$ and $\rho_{2}=\frac{3^{\frac{1}{6}}}{e}\sim \frac{1.2}{e}$. Then 
$$\sqrt{2\pi}\rho_{1}^{k}k^{k}\leq k!\leq \sqrt{2\pi}\rho_{2}^{k}k^{k}$$
for all positive integer $k$.
\end{lem}
{\bf Proof}. Note that 
$$\frac{1}{e},\frac{\sqrt{2}}{e^{2}},\frac{\sqrt{3}}{e^{3}},\frac{\sqrt{4}}{e^{4}},\frac{\sqrt{5}}{e^{5}},\ldots $$
can be expressed by
$$\frac{1}{e},(\frac{2^{\frac{1}{4}}}{e})^{2},(\frac{3^{\frac{1}{6}}}{e})^{3},(\frac{4^{\frac{1}{8}}}{e})^{4},(\frac{5^{\frac{1}{10}}}{e})^{5},\ldots $$
Now,
$$\sup \left\{ \frac{1}{e},\frac{2^{\frac{1}{4}}}{e},\frac{3^{\frac{1}{6}}}{e},\frac{4^{\frac{1}{8}}}{e},\frac{5^{\frac{1}{10}}}{e},\ldots \right\} =\frac{3^{\frac{1}{6}}}{e}$$
implies that $\frac{\sqrt{k}}{e^{k}}\leq \rho_{2}^{k}$ for all $k$. Thus $k!\sim \sqrt{2\pi}\frac{\sqrt{k}}{e^{k}}\cdot k^{k}\leq \sqrt{2\pi}\rho_{2}^{k}\cdot k^{k}$.

The remaining part $\sqrt{2\pi}\rho_{1}^{k}k^{k}\leq k!$ follows by observing that
$$\sqrt{2\pi}(\frac{1}{e})^{k}k^{k}\leq \sqrt{2\pi}(\frac{1}{e})^{k}\cdot \sqrt{k}\cdot k^{k}\sim k!. $$ \hspace{14cm}  $\sharp$
\begin{lem}
  Let $\rho=\frac{\sqrt{3}}{e}$. Then $k!\leq \sqrt{2\pi}\rho^{k}k^{k-1}$ for all $k\geq 1$.
\end{lem}
{\bf Proof}. First,
$$k! \sim \sqrt{2\pi}(\frac{1}{e})^{k}\cdot \sqrt{k}\cdot k^{k}=\sqrt{2\pi}\cdot \frac{k^{\frac{3}{2}}}{e^{k}}\cdot k^{k-1}.$$
Note that $\left\{ \frac{k^{\frac{3}{2}}}{e^{k}}: k=1,2,3,\ldots \right\}$ can be expressed by 
$$\left\{ \frac{1}{e},\frac{2^{\frac{3}{4}}}{e},\frac{3^{\frac{1}{2}}}{e},\frac{4^{\frac{3}{8}}}{e},\ldots \right\}.$$
Our lemma follows by noting that 
$$\sup \left\{ \frac{1}{e},\frac{2^{\frac{3}{4}}}{e},\frac{3^{\frac{1}{2}}}{e},\frac{4^{\frac{3}{8}}}{e},\ldots \right\} =\frac{\sqrt{3}}{e}.$$
\begin{lem}
  Let $h(x)=e^{-\beta|x|^{2}},\beta>0$, be the Gaussian function in $R^{n}$, and $\mu$ be the\\  measure defined in (6). For any positive even integer $l$,
$$\int_{R^{n}}|\xi|^{l}d\mu(\xi)\leq\pi^{\frac{n+1}{2}}\cdot n\cdot \alpha_{n}\cdot 2^{\frac{l+n+2}{2}}\cdot \rho^{\frac{l+n-1}{2}}\cdot \beta^{\frac{l}{2}}\cdot (l+n-1)^{\frac{l+n-3}{2}}\cdot (2+\frac{1}{e})$$
for odd $n$, and
$$\int_{R^{n}}|\xi|^{l}d\mu(\xi)\leq\pi^{\frac{n+1}{2}}\cdot n\cdot \alpha_{n}\cdot 2^{\frac{l+n+3}{2}}\cdot \rho^{\frac{l+n-2}{2}}\cdot \beta^{\frac{l}{2}}\cdot (l+n-2)^{\frac{l+n-4}{2}}$$
for even $n$, where $\rho=\frac{\sqrt{3}}{e}$ and $\alpha_{n}$ is the volume of the unit ball in $R^{n}$.
\end{lem}
{\bf Proof}. 
\begin{eqnarray*}
  &   & \int_{R^{n}}|\xi|^{l}d\mu(\xi) \\
  & = & (\frac{\pi}{\beta})^{\frac{n}{2}}\int_{R^{n}}\frac{|\xi|^{l}}{e^{\frac{|\xi|^{2}}{4\beta}}}d\xi \\
  & = & (\frac{\pi}{\beta})^{\frac{n}{2}}\cdot n\cdot \alpha_{n}\cdot \int_{0}^{\infty}\frac{r^{l}\cdot r^{n-1}}{e^{\frac{r^{2}}{4\beta}}}dr \\
  & = & (\frac{\pi}{\beta})^{\frac{n}{2}}\cdot n\cdot \alpha_{n}\cdot (2\sqrt{\beta})^{l+n}\cdot \int_{0}^{\infty}\frac{r^{l+n-1}}{e^{r^{2}}}dr\\
  & = & \left\{\begin{array}{lll}
                 \pi^{\frac{n}{2}}\cdot n\cdot \alpha_{n}\cdot2^{n-1+l}\cdot \beta^{\frac{l}{2}}\cdot (\frac{l+n-2}{2})(\frac{l+n-2}{2}-1)\cdots (\frac{1}{2})\int_{0}^{\infty}\frac{1}{\sqrt{v}e^{v}}dv \\ 
                 if\  n\  is\  odd, \\ 
                 \pi^{\frac{n}{2}}\cdot n\cdot \alpha_{n}\cdot 2^{n-1+l}\cdot \beta^{\frac{l}{2}}(\frac{l+n-2}{2})!\ if\  n\  is\  even.
               \end{array} \right. 
\end{eqnarray*}
Let $b=\lceil \frac{l+n+2}{2} \rceil$. Then $b! \leq \sqrt{2\pi}\rho^{b}\cdot b^{b-1}$. Then for odd $n$, $b=\frac{l+n+1}{2}$. Thus
$$b! \leq \sqrt{2\pi}\rho^{\frac{l+n-1}{2}}\cdot (\frac{l+n-1}{2})^{\frac{l+n-3}{2}}=\sqrt{2\pi}\cdot 2^{\frac{-l-n+3}{2}}\cdot \rho^{\frac{l+n-1}{2}}\cdot (l+n-1)^{\frac{l+n-3}{2}},$$
and
$$\int_{R^{n}}|\xi |^{l}d\mu(\xi)\leq \pi^{\frac{n+1}{2}}\cdot n\cdot \alpha_{n}\cdot (2+\frac{1}{e})\cdot 2^{\frac{l+n+2}{2}}\cdot \beta^{\frac{l}{2}}\cdot \rho^{\frac{l+n-1}{2}}\cdot (l+n-1)^{\frac{l+n-3}{2}},$$
by noting that $\int_{0}^{\infty}\frac{1}{\sqrt{v}e^{v}}dv\leq 2+\frac{1}{e}$. For even $n$, $b=\frac{l+n-2}{2}$. Thus
$$b!\leq \sqrt{2\pi}\rho^{\frac{l+n-2}{2}}\cdot (\frac{l+n-2}{2})^{\frac{l+n-4}{2}}= \sqrt{2\pi}\cdot 2^{\frac{-l-n+4}{2}}\cdot \rho^{\frac{l+n-2}{2}}\cdot (l+n-2)^{\frac{l+n-4}{2}},$$
and
$$\int_{R^{n}}|\xi|^{l}d\mu(\xi)\leq \pi^{\frac{n+1}{2}}\cdot n\cdot \alpha_{n} \cdot 2^{\frac{l+n+3}{2}}\cdot \beta^{\frac{l}{2}}\cdot \rho^{\frac{l+n-2}{2}}\cdot (l+n-2)^{\frac{l+n-4}{2}}.$$
 \hspace{15cm} $\sharp$

Our interpolation is based on a function space called {\bf native space}, denoted by ${\cal C}_{h,m}$. If
$${\cal D}_{m}=\left\{ \phi \in {\cal D}: \int x^{\alpha}\phi(x)dx=0\ \ for\ \ all\ \ |\alpha|<m\right\}$$, then ${\bf {\cal C}_{h,m}}$ is the class of those continuous functions $f$ which satisfy
\begin{equation}
  \left| \int f(x)\phi(x)dx\right| \leq c(f)\left\{ \int h(x-y)\phi(x)\overline{\phi(y)}dxdy\right\} ^{\frac{1}{2}}
\end{equation}
for some constant $c(f)$ and all $\phi$ in ${\cal D}_{m}$. If $f\in {\cal C}_{h,m}$, let $\| f\| _{h}$ denotes the smallest constant $c(f)$ for which (7) is true. It can be shown that $\| f\| _{h}$ is a semi-norm and ${\cal C}_{h,m}$ is a semi-Hilbert space; in the case $m=0$ it is a norm and a Hilbert space respectively. In this paper both the interpolating and interpolated functions belong to the native space.

The function space ${\cal C}_{h,m}$ is introduced by Madych and Nelson in \cite{MN1} and \cite{MN2}. Later Luh makes a lucid characterization in \cite{Lu1} and \cite{Lu2}. Although there is an equivalent expression for ${\cal C}_{h,m}$ which is easier to understand, we still adopt Madych and Nelson's definition to show the author's respect for them. 
The main result of this paper is the following theorem.
\begin{thm}
  Let $h(x):=e^{-\beta|x|^{2}}$ be the Gaussian function in $R^{n}$. For any positive number $b_{0}$, there are positive constants $\delta_{0},c_{1},c_{2}$, and $c_{3}$ independent of $n$, for which the following is true: If $f\in {\cal C}_{h,m}$ , the native space induced by $h$, and $s$ is the $h$ spline that interpolates $f$ on a subset $X$ of $R^{n}$, then 
\begin{equation}
  |f(x)-s(x)|\leq c_{1}\sqrt{\delta}(c_{2}\delta)^{\frac{c_{3}}{\delta}}\cdot \| f\| _{h}
\end{equation}
for all $x$ in a subset $\Omega$ of $R^{n}$, and $0<\delta \leq \delta_{0}$, where $\Omega$ satisfies the property that for any $x$ in $\Omega$ and any number $\frac{b_{0}}{2}\leq r\leq b_{0}$, there is an $n$ simplex $Q$ with diameter $diamQ=r,x\in Q\subseteq \Omega$, such that for any integer $l$ with $\frac{b_{0}}{\delta}\leq l\leq \frac{2b_{0}}{\delta}$, there is on $Q$ an equally spaced set of centers from $X$ of degree $l-1$.(In fact, the set $X$ can be chosen to consist of these equally spaced centers in $Q$ only.) Here $\|f\| _{h}$ is the $h$-norm of $f$ in the native space. The numbers $\delta_{0},c_{1},c_{2}$, and $c_{3}$ are given by
$ \delta_{0}:=\min \left\{ b_{0},\ \frac{1}{\rho_{3}^{4}\cdot 3^{3}\cdot 2^{7}\cdot b_{0}^{3}}\right\} \ \ where\ \ \rho_{3}=12^{\frac{1}{4}}\cdot \sqrt{e\beta}\ \ $; 
$\left\{
\begin{array}{lll}
 c_{1}:=\left\{ \begin{array}{ll}
                                                                                          \Delta''\cdot \frac{1}{\sqrt{16\pi}}\cdot \frac{1}{\sqrt{b_{0}}}\ \ for\ \ odd\ \ n, \\ \Delta''\cdot \frac{1}{\sqrt{16\pi}}\cdot \frac{1}{\sqrt{b_{0}}}\ \ for\ \ even\ \ n,
                \end{array}\right. \  \\ c_{2}:=\rho_{3}^{4}\cdot 3^{3}\cdot 2^{7}\cdot b_{0}^{3}\ \ , \\ c_{3}:=\frac{b_{0}}{4},
 \end{array} \right.$
 where $\Delta''$ is defined by
$\Delta'':=\left\{ \begin{array}{ll}
                     \sqrt{2+\frac{1}{e}}\cdot \pi^{\frac{n-1}{4}}\cdot (n\alpha_{n})^{\frac{1}{2}}\cdot 2^{\frac{n}{4}}\cdot \rho^{\frac{n-1}{4}}\ for\ odd\ n, \\ \pi^{\frac{n-1}{4}}\cdot (n\alpha_{n})^{\frac{1}{2}}\cdot 2^{\frac{n+1}{4}}\cdot \rho^{\frac{n-2}{4}}\ for\ even\ n,\ with\ \rho=\frac{\sqrt{3}}{e}\ as\ in\ Lemma2.3   
                   \end{array} \right. $ ,where the number $\alpha_{n}$ denotes the volume of the unit ball in $R^{n}$.

In particular, if the point $x$ in $\Omega$ is fixed, the only requirement for $\Omega$ is the existence of an n simplex $Q$, with $diamQ=r,\ x\in Q\subseteq \Omega$, satisfying the afore-mentioned property of equally spaced centers.
\end{thm}
{\bf Proof}. For any $b_{0}>0$, let $\delta_{0}:=\min \left\{ b_{0},\ \frac{1}{\rho_{3}^{4}\cdot 3^{3}\cdot 2^{7}\cdot b_{0}^{3}}\right\} $ where $\rho_{3}=12^{\frac{1}{4}}\cdot \sqrt{e\beta}$. For any $0<\delta\leq \delta_{0}$, there exists an integer $l$ such that $1\leq \frac{\delta}{b_{0}}l\leq 2$ since $0<\frac{\delta}{b_{0}}\leq 1$. Such $l$ satisfies $\frac{b_{0}}{\delta}\leq l\leq \frac{2b_{0}}{\delta}$ and $\frac{b_{0}}{2}\leq \frac{\delta}{2}l\leq b_{0}$.

For any $x\in \Omega$, let $Q$ be an $n$ simplex containing $x$ such that $Q\subseteq \Omega$ and has diameter $diamQ=\frac{\delta l}{2}$. Then Theorem4.2 of \cite{MN2} implies that
\begin{equation}
  |f(x)-s(x)|\leq c_{l}\| f\| _{h}\int_{R^{n}}|y-x|^{l}d|\sigma|(y)
\end{equation}
whenever $l>0$, where $\sigma$ is any measure supported on $X$ such that 
\begin{equation}
  \int_{R^{n}}p(y)d\sigma (y)=p(x)
\end{equation}
for all polynomials $p$ in $P_{l-1}^{n}$. Here
$$c_{l}=\left\{ \int_{R^{n}}\frac{|\xi|^{2l}}{(l!)^{2}}d\mu(\xi)\right\} ^{\frac{1}{2}} $$
whenever $l>0$. Be careful. We temporarily use $c_{1},c_{2},c_{3}$ to denote numbers which are totally different from the $c_{1},c_{2},c_{3}$ mentioned in the theorem.

Let $Y$ be the set of equally spaced centers from $X$ of degree $l-1$ on $Q$. By Lemma1.2, there is a measure $\sigma$ supported on $Y$ such that (10) is satisfied and 
$$\int d|\sigma|(y)\leq \left( \begin{array}{c}
                                 2l-3 \\ l-1    
                               \end{array} \right) $$
The crux of our error estimate is to bound
$$I=c_{l}\int_{R^{n}}|y-x|^{l}d|\sigma|(y).$$
Now Lemma2.3 applies. For odd $n$,
\begin{eqnarray*}
  c_{l} & = & \frac{1}{l!}\left\{ \int_{R^{n}}|\xi|^{2l}d\mu(\xi)\right\} ^{1/2}\\
        & \leq & \frac{1}{l!}\left\{ \pi^{\frac{n+1}{2}}\cdot n\cdot \alpha_{n}\cdot 2^{\frac{2l+n+2}{2}}\cdot \rho ^{\frac{2l+n-1}{2}}\cdot \beta^{l}\cdot (2l+n-1)^{\frac{2l+n-3}{2}}\cdot (2+\frac{1}{e})\right\} ^{1/2}\\
        & \leq & \frac{1}{l!}\left\{ \pi^{\frac{n+1}{4}}\cdot (n\alpha_{n})^{1/2}\cdot 2^{\frac{n+2}{4}}\cdot \rho^{\frac{n-1}{4}}\cdot (\sqrt{2\rho \beta})^{l}\cdot (2l+n-1)^{\frac{2l+n-3}{4}}\cdot \sqrt{2+\frac{1}{e}}\right\} \\
        & \leq & \left\{ \pi^{\frac{n+1}{4}}(n\alpha_{n})^{1/2}\cdot 2^{\frac{n+2}{4}}\cdot \rho^{\frac{n-1}{4}}(\sqrt{2\rho \beta})^{l}(2l+n-1)^{\frac{2l+n-3}{4}}(2+\frac{1}{e})\right\} / \left\{ \sqrt{2\pi}\cdot \rho_{1}^{l}\cdot l^{l}\right\} \\
        &   & by\ Lemma2.1\\
        & = & \pi ^{\frac{n-1}{4}}\cdot (n\alpha_{n})^{1/2}\cdot 2^{\frac{n}{4}}\cdot \rho^{\frac{n-1}{4}}\cdot \rho_{3}^{l}\cdot (\frac{1}{l^{l}})(2l+n-1)^{\frac{2l+n-3}{4}}\cdot (2+\frac{1}{e})^{1/2}\\
        &   & where \ \rho_{3}=\frac{\sqrt{2\rho \beta}}{\rho_{1}} \\
        & = & \Delta''\cdot \rho_{3}^{l}\cdot \frac{1}{l^{l}}(2l+n-1)^{\frac{2l+n-3}{4}}where\ \Delta''=\sqrt{2+\frac{1}{e}}\cdot \pi^{\frac{n-1}{4}}(n\alpha_{n})^{1/2}2^{n/4}\rho^{\frac{n-1}{4}}.       
\end{eqnarray*}
Note that the diameter of $Q$ is $\frac{\delta l}{2}$. This gives for $l\geq n-3$
\begin{eqnarray*}
  I & \leq & \Delta''\rho_{3}^{l}l^{-l}(2l+n-1)^{\frac{2l+n-3}{4}}(diamQ)^{l}\left( \begin{array}{c}
                                                                                       2l-3 \\ l-1
                                                                                    \end{array} \right) \\
    & \leq & \Delta''\rho_{3}^{l}l^{-l}(3l)^{\frac{3l}{4}}(\frac{\delta}{2}l)^{l}\left( \begin{array}{c}
                                                                                           2l-3 \\l-1 
                                                                                        \end{array} \right) \ if\ l\geq n-3 \\
    & = & \Delta''\rho_{3}^{l}(3^{3/4})^{l}(\frac{\delta}{2}l^{3/4})^{l}\left( \begin{array}{c}
                                                                                 2l-3 \\l-1
                                                                               \end{array} \right) \ if\ l\geq n-3\\
    & = & \Delta''(\rho_{3}3^{3/4}\frac{\delta}{2}l^{3/4})^{l}\left( \begin{array}{c}
                                                                        2l-3 \\l-1
                                                                     \end{array}  \right) \\  
    & \leq & \Delta''\left\{ \rho_{3}3^{3/4}\cdot \frac{1}{2}\cdot \delta (\frac{2b_{0}}{\delta})^{3/4}\right\} ^{l}\frac{1}{\sqrt{\pi}}\frac{1}{\sqrt{l-1}}4^{l}\cdot \frac{1}{4}\ by \
    Stirling's\ formula\ and\\
    &   & \frac{b_{0}}{\delta}\leq l\leq \frac{2b_{0}}{\delta}\\
    & = & \Delta''\left\{ \rho_{3}3^{3/4}\cdot 2(2b_{0})^{3/4}\delta^{1/4}\right\} ^{l}\frac{1}{\sqrt{16\pi}} \frac{1}{\sqrt{l-1}} \\
    & \leq & \Delta''\left\{ \rho_{3}3^{3/4}\cdot 2\cdot 2^{3/4}\cdot b_{0}^{3/4}\cdot \delta^{1/4}\right\} ^{\frac{b_{0}}{\delta}}\frac{1}{\sqrt{16\pi}}\cdot \frac{1}{\sqrt{l-1}}\\
    &   & where\ \ \rho_{3}\cdot 3^{3/4}\cdot 2\cdot 2^{3/4}\cdot b_{0}^{3/4}\cdot \delta^{1/4}\leq 1 \ iff\ \delta^{1/4}\leq \frac{1}{\rho_{3}3^{3/4}\cdot 2^{7/4}\cdot b_{0}^{3/4}},\ which\\
    &   & is \ guaranteed\ by\ \delta\leq \delta_{0}:=\min \left\{ b_{0},\ \frac{1}{\rho_{3}^{4}\cdot 3^{3}\cdot 2^{7}\cdot b_{0}^{3}}\right\}\\
    & \leq & \Delta''\left\{ \rho_{3}^{4}\cdot 3^{3}\cdot 2^{7}\cdot b_{0}^{3}\cdot \delta \right\} ^{\frac{b_{0}}{4\delta}}\frac{1}{\sqrt{16\pi}}\cdot \frac{1}{\sqrt{l-1}} \\
    & \leq & \Delta''\cdot \frac{1}{\sqrt{16\pi}}\cdot \left\{ G\delta \right\}^{\frac{g}{\delta}}\cdot \frac{1}{\sqrt{l-1}}\ where\ G=\rho_{3}^{4}3^{3}2^{7}b_{0}^{3}\ and \ g=\frac{b_{0}}{4} \\
    & \sim & \frac{\Delta''}{\sqrt{16\pi}}\left\{ G\delta\right\} ^{\frac{g}{\delta}}\frac{1}{\sqrt{l}}\ (if\ l>>0)\\
    & \leq & \frac{\Delta''}{\sqrt{16\pi}}\sqrt{\frac{1}{b_{0}}}\sqrt{\delta}\left\{ G\delta\right\} ^{\frac{g}{\delta}}\ since\ \frac{b_{0}}{\delta}\leq l\leq \frac{2b_{0}}{\delta}\ 
if\ l\geq n-3
\end{eqnarray*}
For even $n$, if $l\geq n-4$,
\begin{eqnarray*}
  I & = & c_{l}\int_{R^{n}}|y-x|^{l}d|\sigma|(y) \\
    & \leq & \Delta''\rho_{3}^{l}\cdot l^{-l}(2l+n-2)^{\frac{2l+n-4}{4}}(dimQ)^{l}\left( \begin{array}{c}
                                                                                           2l-3 \\ l-1
                                                                                         \end{array} \right) \\
    &   & where\ \Delta''=\pi^{\frac{n-1}{4}}\cdot (n\alpha_{n})^{\frac{1}{2}}\cdot 2^{\frac{n+1}{4}}\cdot \rho^{\frac{n-2}{4}}\ and\ \rho_{3}=\frac{\sqrt{2\rho \beta}}{\rho_{1}}\\
    & = & \Delta''\rho_{3}^{l}\cdot l^{-l}(2l+n-2)^{\frac{2l+n-4}{4}}\left( \frac{\delta l}{2} \right) ^{l}\left( \begin{array}{c}
                     2l-3 \\ l-1 
                   \end{array} \right) \\
    & \leq & \Delta'' \rho_{3}^{l}\cdot (3l)^{\frac{3l}{4}}\cdot \delta^{l}\cdot (\frac{1}{2})^{l}\left( \begin{array}{c}
          2l-3 \\l-1 
        \end{array} \right)\ since\ l\geq n-4\\
    & \leq & \Delta'' (\rho_{3}3^{\frac{3}{4}}l^{\frac{3}{4}}\frac{1}{2}\delta )^{l}\frac{1}{\sqrt{\pi}}\frac{1}{\sqrt{l-1}}4^{l-1}\ by\ Stirling's\ formula \\
    & = & \Delta''\frac{1}{\sqrt{\pi}}\cdot \frac{1}{4}(\rho_{3}3^{\frac{3}{4}}\cdot \frac{1}{2}\cdot \delta l^{\frac{3}{4}}4)^{l}\frac{1}{\sqrt{l-1}} \\
    & \sim & \Delta'' \frac{1}{\sqrt{16\pi}}\cdot \frac{1}{\sqrt{l}}(\rho_{3}\cdot 2\cdot 3^{\frac{3}{4}}\delta l^{\frac{3}{4}})^{l}\ (if\ l>>0)\\
    & \leq & \Delta'' \frac{1}{\sqrt{16\pi}}\cdot \frac{1}{\sqrt{l}}\left\{ \rho_{3}\cdot 2 \cdot 3^{\frac{3}{4}}\cdot \delta\cdot \left( \frac{2b_{0}}{\delta} \right) ^{\frac{3}{4}}\right\} ^{l} \ since\ \frac{b_{0}}{\delta}\leq l \leq \frac{2b_{0}}{\delta}\\
    & = & \Delta''\frac{1}{\sqrt{16\pi}}\cdot \frac{1}{\sqrt{l}}\left\{ \rho_{3}\cdot 2^{\frac{7}{4}}\cdot 3^{\frac{3}{4}}\cdot b_{0}^{\frac{3}{4}}\cdot \delta^{\frac{1}{4}}\right\} ^{l}\\
    & \leq & \Delta''\frac{1}{\sqrt{16\pi}}\cdot \frac{1}{\sqrt{l}}\left\{ \rho_{3}\cdot 2^{\frac{7}{4}}\cdot 3^{\frac{3}{4}}\cdot b_{0}^{\frac{3}{4}}\cdot \delta^{\frac{1}{4}}\right\} ^{\frac{b_{0}}{\delta}}\ where \\
    &   & \rho_{3}\cdot 2^{\frac{7}{4}}\cdot 3^{\frac{3}{4}}\cdot b_{0}^{\frac{3}{4}}\cdot \delta^{\frac{1}{4}}\leq 1\ iff\ \delta\leq \frac{1}{\rho_{3}^{4}\cdot 2^{7}\cdot 3^{3}\cdot b_{0}^{3}},\ which\ is \ guaranteed\ by \\
    &   & \delta \leq \delta_{0}:=\min \left\{ \frac{1}{\rho_{3}^{4}\cdot 2^{7}\cdot 3^{3}\cdot b_{0}^{3}},\ b_{0}\right\} \\
    & \leq & \Delta''\frac{1}{\sqrt{16\pi}}\cdot \frac{1}{\sqrt{l}}\left\{ \rho_{3}^{4}\cdot 2^{7}\cdot 3^{3}\cdot b_{0}^{3}\cdot \delta\right\} ^{\frac{b_{0}}{4\delta}}\\
    & = & \Delta''\frac{1}{\sqrt{16\pi}}\cdot \frac{1}{\sqrt{l}}\left\{ G\delta\right\} ^{\frac{g}{\delta}}\ where\ G=\rho_{3}^{4}\cdot 2^{7}\cdot 3^{3}\cdot b_{0}^{3}\ and\ g=\frac{b_{0}}{4}\\
    & \leq & \Delta''\frac{1}{\sqrt{16\pi}}\cdot \frac{1}{\sqrt{b_{0}}}\cdot \sqrt{\delta}\left\{ G\delta\right\} ^{\frac{g}{\delta}}\ since\ \frac{b_{0}}{\delta}\leq l\leq \frac{2b_{0}}{\delta}.
\end{eqnarray*}
Our error bound (8) then follows immediately by letting $\rho_{3}=\sqrt{e\beta}\cdot 12^{\frac{1}{4}}$. \hspace{4cm}  $\sharp$

In the proof of Theorem2.4 although we require that $l\geq n-3$ for $n$ odd and $l\geq n-4$ for $n$ even, besides $l>>0$, it usually causes no trouble because $\delta\rightarrow 0$ implies $l\rightarrow \infty$ and usually $\delta$ is very small. In order to make the theorem easier to understand, we don't put these requirements into the theorem.\\
\\
{\bf Remark}. In theorem2.4 we avoid using the well-known term ``{\bf fill-distance}'' for scattered data approximation because in our approach the data points are not purely scattered. However the number $\delta$ is in spirit equivalent to the fill-distance in the sense that $\delta\rightarrow 0$ iff the fill-distance $d(Q,Y)\rightarrow 0$ where $Y$ is the equally spaced centers from $X$ of degree $l-1$ on $Q$. Note that $d(Q,Y)\rightarrow 0$ iff $l\rightarrow \infty$, and $l\rightarrow \infty$ iff $\delta \rightarrow 0$.

Furthermore, although the function space ${\cal C}_{h,m}$ is defined in a complicated way, there is a simple expression for it. This can be seen in Wendland's book \cite{We} where Madych and Nelson's native space \cite{MN1} and Wu and Schaback's native space \cite{WS} are unified in an elegant way.\\
\section{Comparison}
The exponential-type error bound for Gaussian interpolation raised by Madych and Nelson in \cite{MN3} is of the form
\begin{equation}
  |f(x)-s(x)|\leq a_{1}(a_{2}\delta)^{\frac{a_{3}}{\delta}}\| f\| _{h},\ as\ \delta \rightarrow 0
\end{equation}
where $a_{1}$ is about the same as $c_{1}$ of (8), $a_{3}=\frac{b_{0}}{8\gamma_{n}}$ and $a_{2}=(3^{3/4}\cdot e\cdot \sqrt{2\rho \beta }\cdot \sqrt{n}\cdot e^{2n\gamma_{n}})^{4}\cdot b_{0}^{3}\cdot \gamma_{n}$ where $\rho =\frac{\sqrt{3}}{e}$, $b_{0}$ is the side length of a cube and $\gamma_{n}$ is defined recursively by
$$\gamma_{1}=2,\ \gamma_{n}=2n(1+\gamma_{n-1})\ if\ n>1.$$
These can be seen in \cite{MN3} and \cite{Lu3}. Note that $\gamma_{n}\rightarrow \infty $ rapidly as $n\rightarrow \infty$. The first few examples are
$$\gamma_{1}=2,\ \gamma_{2}=12,\ \gamma_{3}=78,\ \gamma_{4}=632\ and\ \gamma_{5}=6330.$$
This means that for high dimensions, $a_{2}$ will become extremely large and $a_{3}$ extremely small, making the error bound meaningless.

Our new approach avoids this drawback. Although we require that the centers be equally spaced in the simplex, it causes bo troubles at all both theoretically and practically.

\end{document}